\documentclass[12pt,psfig, reqno]{amsart}
\usepackage{amsmath, amsthm, amsfonts, amssymb, latexsym, graphicx}
\usepackage[latin1]{inputenc}
\pdfoutput=1

\newtheorem{theorem}{Theorem}[section]
\newtheorem{lemma}[theorem]{Lemma}

\theoremstyle{definition}

\newtheorem{example}[theorem]{Example}

\newtheorem{remark}[theorem]{Remark}
\numberwithin{equation}{section}

 \allowdisplaybreaks
\begin{document}

\title{Normality through Sharing of pairs of functions with derivatives}
\author[K. S. Charak]{Kuldeep Singh Charak}
\address{
\begin{tabular}{lll}
& Kuldeep Singh Charak\\
& Department of Mathematics\\
& University of Jammu\\
& Jammu-180 006\\
& India
\end{tabular}}
\email{kscharak7@rediffmail.com }

\author[M. Kumar]{Manish Kumar}
\address{
\begin{tabular}{lll}
& Manish Kumar\\
& Department of Mathematics\\
& University of Jammu\\
& Jammu-180 006\\ 
& India\\
\end{tabular}}
\email{manishbarmaan@gmail.com}

\author[A. Singh]{Anil Singh}
\address{
\begin{tabular}{lll}
& Anil Singh\\
& Department of Mathematics\\
& Maulana Azad Memorial College,\\
& Jammu-180 006\\
& India
\end{tabular}}
\email{anilmanhasfeb90@gmail.com }

\begin{abstract} Let $\mathcal{F}\subset\mathcal{M}(D)$ and let $a, b$ and $c$ be three distinct complex numbers. If, there exist a holomorphic function $h$ on $D$ and a positive constant $\rho$ such that for each $f\in\mathcal{F},$ $f$ and $f^{'}$ partially share three pairs of functions $(a,h), \ (b, c_f)$ and $(c,d_f)$ on $D,$ where $c_f$ and $d_f$ are some values in some punctured disk $D^*_{\rho}(0),$ then $\mathcal{F}$ is normal in $D$. This is an improvement of Schwick's result[Arch. Math. (Basel), \textbf{59} (1992), 50-54]. We also obtain several normality criteria which significantly improve the existing results and examples are given to establish the sharpness of  results.
\end{abstract}

\renewcommand{\thefootnote}{\fnsymbol{footnote}}
\footnotetext{The work of the first author is partially supported by Mathematical Research Impact Centric Support (MATRICS) grant, File No. MTR/2018/000446, by the Science and Engineering Research Board (SERB), Department of Science and Technology (DST), Government of India.}

\maketitle

\section{Introduction } 
Let $D\subseteq \mathbb{C}$ be a domain. For the sake of convenience we shall denote by $\mathcal{M}(D)$ the class of all meromorphic functions on $D$, by $\mathcal{H}(D)$ the class of all holomorphic functions on $D$, and by $\mathbb D$ the open unit disk in $\mathbb C$. Let $f\in \mathcal{M}(D)$ and $a\in \mathbb C.$ Further, we shall denote by $E_f(a)$ the set of $a-$points of $f$. Let $a,b \in \mathbb C$. We say that two functions $f, \ g \in \mathcal{M}(D)$ partially share a pair $(a,b)$ if $z\in E_f(a) \Rightarrow z\in E_g(b)$. Further, if $E_f(a)= E_g(b),$ then $f$ and $g$ are said to share the pair $(a,b)$. Clearly, $f$ and $g$ share the value $a$ if they share the pair $(a,a)$. 

A family $\mathcal F \subset \mathcal{M}(D)$ is said to be normal if each sequence in $\mathcal F$  has a subsequence which converges locally uniformly in $D$ with respect to the spherical metric. The limit function lies in $\mathcal{M}(D)\cup \{\infty\}.$

Mues and Steinmetz \cite{mues} proved that {\it if $f$ is meromorphic in the plane and if $f$ and $f'$ share three values, then $f' \equiv f .$}  Let $\mathcal{F}$ be a subfamily of $\mathcal{M}(D)$ such that for each $f\in \mathcal{F}$, $f$ and $f^{\prime}$ share three distinct values. In view of Bloch's principle a natural question arises: Can $\mathcal{F}$ be normal in $D$? Schwick \cite{schwick1} answered this question affirmatively: 
\begin{theorem}\label{U4:TheoremA} Let $\mathcal F\subset\mathcal{M}(D)$ and let $a, b $ and $c$ be three distinct complex numbers. If, for each $f\in\mathcal F$, $f$ and $f^{'}$ share three pairs of values $(a,a), \ (b,b)$ and $(c,c)$, then $\mathcal F$ is normal in $D$.
\end{theorem} 
 
 Several extensions, improvements and related variants of Theorem \ref{U4:TheoremA} have been obtained by various authors, for example one can see \cite{nevo, li, pang, xu}. The purpose of this paper is to obtain further improvements of results of  Xu \cite{xu} and Li and Yi \cite{li}. 
 \section{Statements of Results }
 Xu \cite{xu} proved that for holomorphic version of Theorem \ref{U4:TheoremA}, the sharing of two distinct values is sufficient to ensure the normality:
\begin{theorem}\label{U4:TheoremB}  Let $\mathcal{F}\subset\mathcal{H}(D)$, and let $a$ and $b$ be two distinct complex numbers. If for each $f\in\mathcal{F}$, $f$ and $f'$ share the pairs of values $(a,a)$ and $(b,b)$, then $\mathcal{F}$ is normal in $D$.
\end{theorem} 
L\"{u}, Xu and Yi \cite{lu} proved Theorem \ref{U4:TheoremB} by  using partial sharing of values: 
\begin{theorem}
\label{U4:TheoremC}  Let $\mathcal{F}\subset\mathcal{H}(D)$, and let $a$ and $b$ be distinct complex numbers. If for each $f\in\mathcal{F}$,  $f$ and $f'$ partially share the pairs of values $(a,a)$ and $(b,b)$, then $\mathcal{F}$ is normal in $D$.
\end{theorem}

Our variation of Theorem \ref{U4:TheoremA} is
\begin{theorem} \label{U4:MTA1}Let $\mathcal{F}\subset\mathcal{M}(D)$ and let $a, b$ and $c$ be three distinct complex numbers. If, there exist a holomorphic function $h$ on $D$ and a positive constant $\rho$ such that for each $f\in\mathcal{F},$ $f$ and $f^{'}$ partially share three pairs of functions $(a,h), \ (b, c_f)$ and $(c,d_f)$ on $D,$ where $c_f$ and $d_f$ are some values in a punctured disk $D^*_{\rho}(0),$ then $\mathcal{F}$ is normal in $D$.
\end{theorem}

\medskip

The values $c_f$ and $d_f$ in Theorem \ref{U4:MTA1} need to be in a finite punctured  disk as shown by the following example:
\begin{example} Consider the family $\mathcal{F}:=\left\{f_n(z)=\tan{nz}:n\in\mathbb{N}\right\}$ of meromorphic functions in $\mathbb{D}.$ Then each $f_n$ and $f'_n$  partially share the pairs $(i,h),$ $(-i,c_f)$ and $(1,2n).$ Note that the values $d_{f_n}=2n$ do not lie in any given finite punctured disk; here $c_{f_n}$ is any value in $\mathbb{C}.$ But $\mathcal{F}$ fails to be normal in $\mathbb{D}.$
\end{example}

\medskip

 The following example shows that the three pairs of functions in Theorem \ref{U4:MTA1}  can not be replaced by two pairs of functions:
\begin{example} Consider the family 
$$\mathcal F:=\left\{f_n(z)=\frac{e^{nz}}{1+e^{nz}}:n\geq 4\right\}\subset\mathcal{M}(\mathbb{D}).$$
 Note that each $f\in\mathcal F$ omits $0$ and $1$ in $\mathbb D$ and for  each $f\in \mathcal F ,$ $f$ and $f'$ partially shares the pairs of functions $(0, h)$ and $(1, c_f),$ where $h$ can be any holomorphic function and $c_f\in\mathbb{C}.$ But the family $\mathcal F$ is not normal in $\mathbb D$.
\end{example}

The holomorphic version of  Theorem \ref{U4:MTA1} is 
\begin{theorem}\label{U4:MTA1holomorphic} Let $\mathcal{F}\subset\mathcal{H}(D)$ and let $a$ and $b$ be two distinct complex numbers. If there exist a holomorphic function $h$ on $D$ and positive constant $\rho$ such that for each $f\in\mathcal{F}$ , $f$ and $f'$ partially share the two pairs  $(a,h)$ and $(b, c_f)$, where $c_f\in D^*_{\rho}(0),$ 
then $\mathcal{F}$ is normal in $D$.
\end{theorem}
 Note that Theorem \ref{U4:MTA1holomorphic} is an improvement of Theorem \ref{U4:TheoremC}. The values $ c_f$ in Theorem \ref{U4:MTA1holomorphic}  have to be essentially in a finite punctured disk, which is clear from the following example:
\begin{example} Consider the family 
$$\mathcal{F}:=\left\{f_n=e^{nz}:n\in\mathbb{N}\right\}\subset\mathcal{H}(\mathbb{D}).$$
  Then $f_n$ and $f_n'$ partially share the pairs $(0,0)$ and $(1,n)$. Note that $c_{f_n}=n$ are not contained in any finite disk and the family $\mathcal F$ is not normal in $\mathbb{D}$.
\end{example}

Li and  Yi \cite{li} considered partial sharing of the pair of values $(a, a)$ by $f $ and $f'$ and another  pair of values $(b,b)$  partially shared by $f'$ and $f,$ and obtained the following normality criterion:
  
\begin{theorem}\label{U4:TheoremD}
Let $\mathcal{F}\subset\mathcal{H}(D)$  and let  $a,$  $b\in\mathbb{C}$ be distinct  such that $b\neq 0$. If for each $f\in\mathcal F,$   $f$ and $f'$ partially share the pair $(a,a)$  and   $f'$ and $f$ partially share the pair $(b,b)$, then $\mathcal F$ is normal in $D$.
  \end{theorem}

 Let $A\subset {D}$ and $a\in \mathbb{C}.$ For $f,g\in\mathcal{M}(D),$ we shall say that $f$ and $g$ partially share the pair $(a, A),$ if $f(z)=a$ implies $g(z)\in A.$
 
 As an improvement of Theorem \ref{U4:TheoremD}, we have obtained the following result:
\begin{theorem} \label{U4:MTA2} Let $\mathcal{F}\subset\mathcal{H}(D)$, and let  $a$ and $b\neq 0$ be two distinct complex numbers. Let $A$ be a compact set such that $b\notin A$ and $B=\left\{z:|z-a|\geq\epsilon\right\}$, for some $\epsilon>0$. If for each $f\in \mathcal{F},$ $f$ and $f'$ partially share the pair $(a, A)$ and  $f'$ and $f$ partially share the pair $(b, B)$, then $\mathcal{F}$ is normal in $D$.
\end{theorem}

\begin{remark}
After obtaining Theorem \ref{U4:MTA2} as an improvement of Theorem \ref{U4:TheoremD} we came across a result of  Sauer and Schweizer \cite{Andreas-1}: {\it Let $\mathcal F$ be a family of holomorphic functions in a domain $D.$ Let $a$ and $b\neq 0$ be two complex numbers such that $b\neq a$, and let $A$ and $B$ be compact subsets of $\mathbb C$ with $b\notin A$ and $a\notin B$. If, for each $f\in \mathcal F$ and $z\in D$, $f$ and $f^{'}$ partially share the pair $(a, A)$ and $f^{'}$ and $f$ partially share the pair $(b, B)$, then $\mathcal F$ is normal in $D$.} This result is  also an improvement of Theorem \ref{U4:TheoremD}. Theorem \ref{U4:MTA2} also provides an improvement of Sauer and Schweizer's result.
\end{remark}

The condition `the set $B$ must be at a positive distance away from the point $a$' in Theorem \ref{U4:MTA2} cannot be dropped as shown by the following example:

\begin{example}
Let $\mathcal F:=\left\{e^{nz}:n\in\mathbb N\right\}\subset \mathcal H(\mathbb D)$. Take $a=0$ and $b=1$. Then $f(z)\neq a$ and $f^{'}_n(z)=b\Rightarrow f_n(z)=1/n\to a$. But $\mathcal F$ is not normal at $z=0$.
\end{example}

In the next example, we show that the boundedness of set $A$ in Theorem \ref{U4:MTA2} can not be relaxed:

\begin{example}
Let $\mathcal F:=\left\{e^{nz}/n:n\in\mathbb N\right\}\subset \mathcal H(\mathbb D)$. Take $a=1$ and $b=-1$. Then $f_n(z)=1\Rightarrow f^{'}_n(z)=n\in\mathbb N$ and $f^{'}_n(z)=-1\Rightarrow f_n(z)=-1/n\in\{z:|z-1|\geq 1\}$. But $\mathcal F$ is not normal at $z=0$.
\end{example}

Another variant of Theorem \ref{U4:MTA2} is obtained as: 
\begin{theorem}\label{U4:MTA3}  Let $\mathcal{F}\subset\mathcal{H}(D)$ be such that zeros of each $f\in\mathcal{F}$ have multiplicity at least $k,$ where $k\in\mathbb{N}$ and $b(\neq 0)\in\mathbb{C}.$ Let $A$ be a compact set and $B=\left\{z:|z|\geq\epsilon\right\}$ for some $\epsilon>0 .$ If for each $f\in\mathcal F,$ $f$ and $f^{(k)}$ partially share the pair $(0, A)$ and $f^{(k)}$ and $f$ partially share the pair $(b, B)$  in $D ,$ then $\mathcal{F}$ is normal in $D.$
\end{theorem}

The condition $`{b\neq 0}$' in Theorem \ref{U4:MTA3} can not be dropped, as can be seen from the following example:
\begin{example} Let 
$\mathcal F:=\left\{e^{nz}:n\in\mathbb N\right\}\subset\mathcal{H}(\mathbb{D}).$ Then $\mathcal F$ satisfies all the conditions of Theorem \ref{U4:MTA3} with $b=0, $ but $\mathcal F$ is not normal in $\mathbb D$.
\end{example}

Also, the condition `the zeros of $f\in\mathcal F$ have multiplicity at least $k$' in Theorem \ref{U4:MTA3} can not be weakened:

\begin{example} Consider the family
 $\mathcal F:=\{n\sinh{z}:n\in \mathbb N\}\subset \mathcal H(\mathbb D).$  Then, clearly, the zeros of $f\in\mathcal F$ are simple and $f\equiv f^{''}$. But the family $\mathcal F$ is not normal at $z=0$.
\end{example}
The meromorphic version of Theorem \ref{U4:MTA2} does not hold as shown by the following example :

\begin{example} Let $a\in\mathbb{C}\setminus \{1\}$ and consider the family
$$\mathcal F:=\left\{\frac{n+(nz-1)^{2}}{n(nz-1)}+a:n\in\mathbb N\right\}\subset\mathcal{M}(\mathbb{D}).$$ 
 Then, for each $f\in\mathcal F$, $f$ and $f'$ partially share the pair $(a,2)$ and  $f^{\prime}\neq 1.$ Thus, for each $f\in \mathcal F,$ $f$ and $f'$ partially shares the pair $(a, A)$ and $f'$ and $f$ partially shares the pair $(1, B), $ where $A=\{2\}$ and $B=\{z:|z-a|\geq\epsilon\}$ for any $\epsilon>0.$
But $\mathcal F$ is not normal in $\mathbb D$.
\end{example}

However, the following related meromorphic version  holds:
\begin{theorem}\label{U4:MTA4} Let $\mathcal{F}\subset\mathcal{M}(D)$ be such that zeros of each $f\in\mathcal{F}$ have multiplicity at least $k+1,$ where $k\in\mathbb{N} .$  Let $a$ and $b$ be two distinct non-zero complex numbers, and $A$ be a compact set and $B=\left\{z\in \mathbb{C}:|z|\geq \epsilon\right\}$ for some $\epsilon>0$. If for each $f\in\mathcal F$, $f$ and $f^{(k)}$ partially share the pair $(a, A)$ and $f^{(k)}$ and $f$ partially share the pair $(b, B)$, then $\mathcal F $ is normal in $D$.
\end{theorem}


The following example shows that the condition `zeros of each $f\in\mathcal{F}$ have multiplicity at least $k+1,$' in Theorem \ref{U4:MTA4} is essential:

\begin{example}
Consider the family
$$\mathcal F:=\left\{\frac{e^{nz}}{n}+2:n\in\mathbb N\right\}.$$
of entire functions. Then, clearly, $f(z)\neq 2$ and $f^{'}(z)=1\Rightarrow f(z)=1/n+2\in \{z:|z|\geq 2\}$. Since $f^{'}(z)\neq 0$, the zeros of $f$ are simple. But the family $\mathcal F$ is not normal at $z=0$. 
\end{example}

Also, the condition `set $B$ must be at a positive distance away from the origin' in Theorem \ref{U4:MTA4} cannot be dropped:
\begin{example}
Consider the family
$$\mathcal F:=\left\{\frac{1}{e^{nz}+1}:n\in\mathbb N\right\}\subset \mathcal M(\mathbb D).$$
Take $a=1,\ b=-1$. Then, clearly, $f\neq0, 1$. Also, $$f^{'}_n(z)=-1\Rightarrow f_n(z)=\frac{2}{\left\{(n-2)\pm\sqrt{(n-2)^2-4}\right\}+2}$$
which are not contained in any set of the form $\{z:|z|\geq\epsilon\}$, for any $\epsilon>0$. But the family $\mathcal F$ is not normal at $z=0$. 
\end{example}

 \section{Proofs of the results}

To prove the results of this paper, we require the following  lemmas:
\begin{lemma} \cite{pang}\label{U4:l1} Let $\mathcal{F}\subset\mathcal M(\mathbb D)$ be such that for each $f\in \mathcal F$, all zeros of $f$ are of multiplicity at least $k$. Suppose that there exists a number $L\geq 1$ such that $|f^{(k)}(z)|\leq L$ whenever  $f\in \mathcal{F}$ and $f(z)=0$. If $\mathcal{F}$ is not normal in $\mathbb{D}$, then for every $\alpha\in [0, k]$, there exist  $r\in (0,1)$,  $\{z_n\}\subset D_r(0),$ $\{f_n\}\subset \mathcal{F}$ and 
$\{\rho_n\}\subset (0,1):\rho_n\rightarrow 0$ such that
$$g_n(\zeta)=\rho_{n}^{-\alpha}f_n(z_n+\rho_n\zeta) \to g(\zeta)$$
locally uniformly on $\mathbb{C}$ with respect to the spherical metric, where $g$ is a non-constant meromorphic function on $\mathbb C$ with
	$g^{\#}(\zeta)\leq g^{\#}(0)=kL+1.$
\end{lemma}

\begin{lemma} \label{U4:l4}\cite{W3} Let $g\in \mathcal{M}(\mathbb{C})$ be of finite order. If $g$ has only finitely many critical values, then it has only finitely many asymptotic values.
\end{lemma}

\begin{lemma} \label{U4:l5} \cite{W1} Let $g\in \mathcal{M}(\mathbb{C})$ be  transcendental having no poles at the origin and let  the set of finite critical and asymptotic values of $g$ be bounded. Then there exists $R>0$ such that 
$$|g^{\prime}(z)|\geq \frac{|g(z)|}{2\pi|z|}\log{\frac{g(z)}{R}},$$
 for all $z\in\mathbb{C}\setminus \left\{0\right\}$ which are not poles of $g$.
\end{lemma}

\begin{lemma} \label{U4:l6} \cite{W3} Let $f\in\mathcal{M}(\mathbb{C})$ be transcendental and of finite order. Suppose all  zeros of $f$ have multiplicity at least $k+1$, where $k\in\mathbb{N}.$ Then $f^{(k)}$ assumes every non-zero complex number infinitely often.
\end{lemma}
 
\textbf{Proof of Theorem \ref{U4:MTA1}:}
Suppose that $\mathcal{F}$ is not normal. Then $\mathcal{F}_a=\left\{f-a:f\in\mathcal{F}\right\}$ is not normal and therefore, by Zalcman Lemma, there exist a sequence $\left\{f_{n}-a\right\}\subset\mathcal{F}_a$, sequence  $\left\{z_n\right\}$ of points in  $D$ and a sequence $\{\rho_n\}$ of positive real numbers with  $\rho_n\to 0$ as $n\to\infty$ such that the re-scaled sequence $ \{g_n(\zeta):=f_n(z_n+\rho_n\zeta)-a\} $ converges locally uniformly to a  non-constant meromorphic function $g$ on $\mathbb{C} .$ 
 
\medskip

 Suppose $g(\zeta_0)=0$. Then by Hurwitz's Theorem, there exists a sequence $\zeta_n\rightarrow \zeta_0$ as $n\to\infty$ such that for sufficiently large $n,$ $ g_n(\zeta_n)=0.$ That is, 
 $f_n(z_n+\rho_n\zeta_n)=a .$  
Thus, by hypothesis, $f_{n}^{\prime}(z_n+\rho_n\zeta_n)=h(z_n+\rho_n\zeta_n)$, and hence
 $$g^{\prime}(\zeta_0)=\lim\limits_{n\to\infty}g^{\prime}_{n}(\zeta_n)=\lim\limits_{n\to\infty}\rho_n f^{\prime}_{n}(z_n+\rho_n\zeta_n)= \lim\limits_{n\to\infty}\rho_n h(z_n+\rho_n\zeta_n)=0.$$ 
 This shows that the zeros  of $g$ have multiplicity at least $2.$  Similarly, we can show that the zeros of $g-(b-a)$ and  $g-(c-a)$ have multiplicity at least $2.$

\medskip 

 Next, we show that $g$ omits $b-a$. Suppose that $\zeta_0$ is a zero of $g-(b-a)$ with multiplicity $k $. Then 
 \begin{equation}\label{eq:3.1}
 g^{(k)}(\zeta_0)\neq 0. 
 \end{equation}
 
 Choose $\delta>0$ such that 
\begin{equation}\label{eq:4}
g(\zeta)\neq b-a,\  g^{\prime}(\zeta)\neq 0, \cdots,\ g^{(k)}(\zeta)\neq 0
\end{equation}
on $D^{*}_{\delta}(\zeta_0).$ 

Since $g(\zeta_0)=b-a$, by Hurwitz's Theorem, there exists $\zeta_{n,i}\to\zeta_0, n\to\infty\ (i=1,\cdots ,k)$ in $D_{\delta}(\zeta_0)$ such that $g_{n}(\zeta_{n,i})=b-a ,$ for sufficiently large $n .$ That is, $f_{n}(z_n+\rho_n\zeta_{n,i})=b$ and thus $0<|f'_{n}(z_n+\rho_n\zeta_{n,i})|\leq \rho .$

Further,
\begin{equation}\label{eq:5}
g^{\prime}_{n}(\zeta_{n,i})=\rho_nf^{\prime}_{n}(z_n+\rho_n\zeta_{n,i})\neq 0, \mbox{ for }i=1,\cdots ,k .
\end{equation}
This implies $\zeta_{n,i}, (i=1,2,\cdots,k)$ are  simple zeros of $g_{n}-(b-a).$
 
 Also $\zeta_{n,i}\neq \zeta_{n,j} \ (1\leq i<j\leq k)$ and 
$$g^{\prime}(\zeta_0)=\lim\limits_{n\to\infty}g^{\prime}_{n}(\zeta_{n,i})=0.$$ 
Therefore, by (\ref{eq:5}), for sufficiently large $n,$ 
$\ g^{\prime}_{n}-\rho_nc_{f_{n}},$ where $c_{f_n}=f_n'(z_n+\rho_n\zeta_{n,i}),$ has at least $k$ zeros $\zeta_{n,i} (i=1,\cdots,k)$ in $D^{*}_{\delta}(0).$  This implies that $\zeta_0$ is a zero of $g^{\prime}$ with multiplicity at least $k$ and hence $g^{(k)}(\zeta_0)=0$, which contradicts (\ref{eq:3.1}). Hence $g(\zeta)\neq b-a$.  Similarly, we can show that $g$ omits $c-a$ and then by second fundamental theorem of Nevanlinna, we arrive at a contradiction. 
$\hfill \Box$
 
 \bigskip

The Proof of Theorem \ref{U4:MTA1holomorphic} is obtained exactly on the lines of the proof of Theorem \ref{U4:MTA1}, so we omit it.

\bigskip

\textbf{Proof of Theorem \ref{U4:MTA2}:}
 We may assume that $D$ is the open unit disk $\mathbb{D}$. Suppose that $\mathcal{F}$ is not normal in $\mathbb{D}$. Then  $\mathcal{F}_a=\left\{f-a:f\in\mathcal{F}\right\}$ is not normal in $\mathbb{D}$. For any $h\in\mathcal{F}_a$,  $|h^{\prime}(z)|\leq M+1$ whenever $h(z)=0$, where $M=\sup \left\{|z|:z\in A\right\}$.  By Lemma \ref{U4:l1}, there exist a sequence $\left\{f_{n}-a\right\}\subset\mathcal{F}_a$, sequence  $\left\{z_n\right\}$ of points in  $D$ and a sequence $\{\rho_n\}$ of positive real numbers with  $\rho_n\to 0$ as $n\to\infty$ such that 
\begin{equation} \label{eq:3}
	g_{n}(\zeta)=\rho_{n}^{-1}\left(f_n(z_n+\rho_n\zeta)-a\right)\rightarrow g(\zeta)
\end{equation}
as $n\to\infty$, locally uniformly on $\mathbb C$, where $g$ is a non-constant entire function satisfying 
$$g^{\#}(\zeta)\leq g^{\#}(0)=M+2$$
implying that the order of $g$ is at most 1.\\

 \textbf{Assertion 1:} If $g(z)=0$, then $g^{\prime}(z)\in A.$\\
Suppose that $g(\zeta_0)=0$. Then by Hurwitz's Theorem, there exists $\zeta_n\to\zeta_0$ as $n\to\infty$ such that for sufficiently large $n,$ $g_n(\zeta_n)=0$ . This implies that $f_n(z_n+\rho_n\zeta_n)=a$. Since $f$ and  $f^{\prime}$ partially share the pair $(a, A)$, 
$$g_{n}^{\prime}(\zeta_n)=f^{\prime}_n(z_n+\rho_n\zeta_n)\in A.$$
Since $A$ is compact, $$g^{\prime}(\zeta_0)=\lim\limits_{n\to \infty}g^{\prime}_n(\zeta_n)\in A$$ 
and this proves Assertion 1. \\

 \textbf{Assertion 2:} $g^{\prime}(\zeta)\neq b,  \forall \ \zeta \in \mathbb{C}.$\\
Suppose that $g^{\prime}(\zeta_0)=b$ for some $\zeta_0\in\mathbb{C} .$  If $g^{\prime}(\zeta)\equiv b$, then $g(\zeta)=b\zeta+c$, so by Assertion 1, $b\in A$, a contradiction. Thus $g^{\prime}(\zeta)\not\equiv b$.

 Now by Hurwitz's Theorem, there exists $\zeta_n\rightarrow \zeta_0$ as $n\to \infty$, such that for sufficiently large $n ,$
$$g^{\prime}_n(\zeta_n)=f^{\prime}_n(z_n+\rho_n\zeta_n)=b.$$
Since  $f^{'}$ and $f$ partially share the pair $(b, B)$, 

$$ |g_n(\zeta_n)|=\rho^{-1}_n|\left(f_n(z_n+\rho_n\zeta_n)-a\right)|\geq \frac{\epsilon}{\rho_n}\to\infty\mbox{ as }n\to\infty.$$
That is, $g(\zeta_0)=\infty$, a contradiction since $g^{\prime}(\zeta_0)=b$. This proves Assertion 2.\\

 Since  $g$ is of order at most $1,$ so is $g^{\prime}$ and then by Assertion 2, we have 
$$g^{\prime}(\zeta)=b+e^{l+m\zeta}.$$
where $l,m\in \mathbb{C}$.

 Now we have the following two cases:

\noindent\textbf{Case-1}. When $m\neq 0.$ In this case, $g$ is a transcendental entire function of order one. Since $g^{\prime}$ omits $b(\neq 0)$,  by Hayman's alternative  $g$ has infinitely many zeros $\left\{z_i\right\}:\ |z_i|\to \infty$ as $i\to \infty$.

 Define $G(z)=g(z)-bz$, then $G^{'}(z)=g^{\prime}(z)-b\neq 0$, $G$ has no critical values. Thus by Lemma \ref{U4:l4}, $G$ has only finitely many asymptotic values. Applying Lemma \ref{U4:l5} to $G$, we have 

$$\frac{|z_iG^{'}(z_i)|}{|G(z_i)|}\geq\frac{1}{2\pi}\log{\frac{|G(z_i)|}{R}}=\frac{1}{2\pi}\log{\frac{|bz_i|}{R}}.$$
This implies

\begin{equation}\label{eq:U4.00}
\frac{|z_iG^{'}(z_i)|}{|G(z_i)|}\to\infty\mbox{ as }i\to\infty.
\end{equation}
 Since $g=0\Rightarrow |g^{\prime}|\leq M ,$ which further implies that ${|z_iG^{'}(z_i)|}/{|G(z_i)|}$ is bounded. Thus (\ref{eq:U4.00}) yields a contradiction.\\

\noindent\textbf{Case-2}. When $m=0.$ In this case $g(\zeta)=\left(b+e^{l}\right)\zeta+t,$ where $t$ is a constant. By Assertion 1, we get $b+e^{l}\in A$. Thus $g^{\#}(0)<M+2$, a contradiction.  $\hfill \Box$
  
\bigskip
 
\textbf{Proof of Theorem \ref{U4:MTA3}:}
We may assume that $D$ is the open unit disk $\mathbb{D}$. Suppose that $\mathcal F$ is not normal in $\mathbb{D}$. Then, by Lemma \ref{U4:l1}, (with $\alpha=k$ and $L=M+1$, where $M=\sup\{|z|:z\in A\}$), there exist $f_n\in\mathcal F,\ z_n\in \mathbb{D}$ and $\rho_n\to 0^+$ such that 
$$g_n(\zeta)=\frac{f_{n}(z_n+\rho_n\zeta)}{\rho_{n}^{k}}\to g(\zeta)$$
locally uniformly on $\mathbb C$, where $g$ is a non-constant entire function such that $g^{\#}(\zeta)\leq g^{\#}(0)=k(M+1)+1$ and the order of $g$ is at most one.

\smallskip

 Next we show that  zeros of $g$ are of multiplicity at least $k$ and $ g(z)=0$ implies that $g^{(k)}(z)\in A .$  Let $g(\zeta_0)=0$. Then by Hurwitz's Theorem, there exists a sequence $\zeta_n\to\zeta_0$ as $n\to\infty$ such that for sufficiently large $n,$ $\ g_n(\zeta_n)=0$ . That is $f_n(z_n+\rho_n\zeta_n)=0$ and by assumption, we have,  $f_{n}^{(i)}(z_n+\rho_n\zeta_n)=0\ (i=1,\cdots, k-1)$ and $f_{n}^{(k)}(z_n+\rho_n\zeta_n)\in A.$
Thus
$$g^{(i)}(\zeta_0)=\lim_{n\to\infty}g_{n}^{(i)}(\zeta_n)=\lim_{n\to\infty}\rho^{i-k}_{n}f_{n}^{(i)}(z_n+\rho_n\zeta_n)=0\ (i=1,\cdots,k-1)$$
and
$$g^{(k)}(\zeta_0)=\lim_{n\to\infty}g_{n}^{(k)}(\zeta_n)=\lim_{n\to\infty}f_{n}^{(k)}(z_n+\rho_n\zeta_n)\in A.$$
Therefore, all zeros of $g$ are of multiplicity at least $k$ and $g(z)=0$ implies that $g^{(k)}(z)\in A$.\\

\noindent\textbf{Assertion:} $g^{(k)}(z)\neq b$ in $\mathbb C$.\\
Suppose that $g^{(k)}(\zeta_0)=b$. If $g^{(k)}(\zeta)\equiv b$, then $g$ is a polynomial of degree $k$. Since all zeros of $g$ are of multiplicity at least $k$, $g$ has only one zero, say $\zeta^{\prime}$. Thus 

$$g(\zeta)=\frac{b(\zeta-\zeta^{\prime})^k}{k!}.$$
Since $g(\zeta)=0 \Rightarrow g^{(k)}(\zeta)\in A$, $|b|\leq M$. By a simple calculation, we have 
$$
 g^{\#}(0)\leq
\begin{cases}
k/2 &; |\zeta^{\prime}|\geq 1\\
M  & ; |\zeta^{\prime}|< 1\

\end{cases}
$$
That is, $g^{\#}(0)<k(M+1)+1$, a contradiction. Thus $g^{(k)}(\zeta)\not\equiv b$.

 Thus, we choose a sequence $\zeta_n\to\zeta_0$ as $n\to\infty$ such that $g^{(k)}_{n}(\zeta_n)=b$. This implies that $f^{(k)}_{n}(z_n+\rho_n\zeta_n)=b$ and by hypothesis, we find that $|f_{n}(z_n+\rho_n\zeta_n)|\geq \epsilon.$
 
Therefore ,
$$\left|g(\zeta_0)\right|=\lim_{n\to\infty}|g_{n}(\zeta_n)|=\lim_{n\to\infty}\left|\frac{f_n(z_n+\rho_n\zeta_n)}{\rho_{n}^{k}}\right|\geq \lim_{n\to\infty}\frac{\epsilon}{\rho^{k}_{n}}= \infty.$$
That is, $g(\zeta_0)=\infty$, a contradiction since $g^{(k)}(\zeta_0)=b$ and this proves the Assertion.\\
Since  $g$ is of order at most one, so is $g^{(k)}$ and by Assertion, we find that 
$$g^{(k)}(\zeta)=b+e^{l+m\zeta},$$
where $l$ and $m$ are constants. Now we have the following two cases:\\

\noindent\textbf{Case-I.} If $m=0$, then $g$ is a polynomial of degree $k$. Since all zeros of $g$ are of multiplicity at least $k$, $g$ has only one zero, say $\zeta^{\prime}$. Thus $$g(\zeta)=\frac{(b+e^{l})(\zeta-\zeta^{\prime})^{k}}{k!}.$$ By second part of Assertion, we have $|b+e^{l}|\leq M$ and as obtained above, we have that $g^{\#}(0)<k(M+1)+1$, a contradiction. \\

\noindent\textbf{Case-II.} If $m\neq 0$. then $g$ is a transcendental entire function. Since $g^{(k)}(\zeta)\neq b(\neq 0)$, by Hayman's alternative, $g$ has infinitely many zeros $\left\{z_i\right\}$ and $|z_i|\to\infty$ as $n\to\infty$. Define $G(z)=g^{(k-1)}(z)-bz$, then $G^{'}(z)=g^{(k)}(z)-b\neq 0$, $G$ has no critical value. Thus by Lemma \ref{U4:l4}, $G$ has only finitely many asymptotic values. Applying Lemma \ref{U4:l5} to $G$, we have

$$\frac{|z_iG^{'}(z_i)|}{|G(z_i)|}\geq\frac{1}{2\pi}\log{\frac{|G(z_i)|}{R}}=\frac{1}{2\pi}\log{\frac{|bz_i|}{R}}.$$
This implies that 

$$\frac{|z_iG^{'}(z_i)|}{|G(z_i)|}\to\infty$$
as $i\to\infty,$ which leads to a contradiction, since $g=0$ implies $g^{(k)}\in A $ and ${|z_iG^{'}(z_i)|}/{|G(z_i)|}$ is bounded.    $\hfill\Box$

\bigskip

\textbf{Proof of Theorem \ref{U4:MTA4}:}
We may take  $D$ to be $\mathbb{D},$ the open unit disk. Suppose that $\mathcal F$ is not normal on $\mathbb{D}$. Then, by Lemma \ref{U4:l1}, there exist $z_n\in \mathbb{D}, f_n\in\mathcal F$ and $\rho_n\to 0^{+}$ such that $\left\{g_n(\zeta)=\rho_n^{-k}\left(f_n(z_n+\rho_n\zeta)\right)\right\}$ converges spherically locally uniformly on $\mathbb{C}$ to a non-constant meromorphic function $g$, all of whose zeros have multiplicity at least $k+1$ and the order of $g$ is finite.

\noindent\textbf{Assertion 1:} $g^{(k)}\neq b$ on $\mathbb C$.

Suppose that $g^{(k)}(\zeta_0)=b$, for some $\zeta_0\in\mathbb C$. If $g^{(k)}\equiv b$, then $g$ is a polynomial of degree $k$, a contradiction since all zeros of $g$ are of multiplicity at least $k+1$. Thus by Hurwitz's Theorem, there exists $\zeta_n\to\zeta_0$ such that for sufficiently large $n$ ,
$$g^{(k)}_n(\zeta_n)=f^{(k)}_n(z_n+\rho_n\zeta_n)=b.$$
By assumption, $|f_n(z_n+\rho_n\zeta_n)|\geq\epsilon$ and so 
$$|g(\zeta_0)|=\lim\limits_{n\to\infty}|g_n(\zeta_n)|=\lim\limits_{n\to\infty}\frac{|f_n(z_n+\rho_n\zeta_n)|}{\rho_n^{k}}\geq\lim\limits_{n\to\infty}\frac{\epsilon}{\rho_n^{k}}=\infty.$$
That is, $g(\zeta_0)=\infty$, a contradiction since $g^{(k)}(\zeta_0)=b$.\\

\textbf{Assertion 2:} $g$ is an entire function.

 Suppose that $g(\zeta_1)=\infty$, for some $\zeta_1\in\mathbb C$. For sufficiently large $n$, we can choose a closed disk $\overline{D_r(\zeta_1)}$ such that $g_n(\zeta)\neq 0$ and $g(\zeta)\neq 0$, and $1/g_n(\zeta)\to 1/g(\zeta) $ uniformly on $\overline{D_r(\zeta_1)}$. Thus 
$$\frac{1}{g_n(\zeta)}-\frac{\rho_n^{k}}{a}\to\frac{1}{g(\zeta)}, $$
uniformly on $\overline{D_r(\zeta_1)}$. Since $1/g(\zeta_1)=0$, there exits $\zeta_n\to\zeta_1$ such that for sufficiently large $n$,
$$\frac{1}{g_n(\zeta_n)}-\frac{\rho_n^{k}}{a}=0.$$ 
That is, $f_n(z_n+\rho_n\zeta_n)=a$. By assumption, we have $|f^{(k)}_n(z_n+\rho_n\zeta_n)|\leq M$, where $M=\sup\{|z|:z\in A\}$ and hence $|g^{(k)}(\zeta_1)|\leq M$, a contradiction since $g(\zeta_1)=\infty$.

 Since $g$ is entire and $g^{(k)}\neq b$ on $\mathbb C$, by Lemma \ref{U4:l6}, $g$ is  a polynomial of degree at most $k$, a contradiction.     $\hfill \Box$

\bigskip
 
\bibliographystyle{amsplain}

\begin{thebibliography}{99}
\bibitem{W1} W. Bergweiler, {\em On the zeros of certain homogeneous differential polynomial}, Arch. Math., \textbf{64} (1995), 199-202.
\bibitem{W3} W. Bergweiler and A. Eremenko, {\em On the singularities of the inverse to a meromorphic function of finite order}, Rev. Mat. Iberoamericana., \textbf{11}  (1995), 355-373.
\bibitem{nevo} J. Grahl and Shahar Nevo, {\em On the result of Singh and Singh concerning shared values and normality}, Complex Variables and Elliptic Equations, \textbf{55(2)} (2010), 347-356.
\bibitem{li} J. Li and H. Yi, {\em Normal families and uniqueness of entire functions and their derivatives}, Arch. Math., \textbf{87} (2006), 52-59.
\bibitem{lu} F. L\"{u}, J. Xu and H. Yi, {\em Uniqueness theorems and normal families of entire functions and their derivatives}, Ann. Pol. Math., \textbf{95}\textbf{(1)} (2009), 67-75.
\bibitem{mues} E. Mues and N. Steinmetz, {\em Meromorphe Funktionen, die mit ihrer AbleitungWerte teilen,}Manuscripta Math., \textbf{29} (1979), 195-206.
\bibitem{pang} X. C. Pang and L. Zalcman, {\em Normal families and shared values}, Bull. London Math. Soc., \textbf{32} (2000), 325-331.
\bibitem{schwick1}W. Schwick, {\em Sharing values and normality,} Arch. Math. (Basel), \textbf{59} (1992), 50-54.

\bibitem{Andreas-1}A. Sauer and A. Schweizer, {\em A uniqueness problem concerning entire functions and their derivatives} arXiv:2208.11341 (2022).

\bibitem{xu} Y. Xu, {\em Normality criteria concerning sharing values}, Indian J. Pure Appl. Math., \textbf{30} (1999), 287-293.






\end{thebibliography}

\end{document}